\documentclass{article}

\usepackage{arxiv}

\usepackage[utf8]{inputenc} 
\usepackage[T1]{fontenc}    
\usepackage{hyperref}       
\usepackage{url}            
\usepackage{booktabs}       
\usepackage{amsfonts}       
\usepackage{nicefrac}       
\usepackage{microtype}      
\usepackage{lipsum}		
\usepackage{graphicx}
\usepackage{natbib}
\usepackage{doi}

\usepackage{amsmath}
\usepackage{amsthm}

\usepackage{amssymb}
\usepackage{algorithm}
\usepackage{algpseudocode}
\usepackage{tikz}
\usepackage{graphicx,subfigure}
\usepackage{makecell}
\usetikzlibrary{shapes.arrows}
\usetikzlibrary{calc}
\usepackage{array,multirow}
\DeclareMathOperator*{\minimize}{minimize}
\DeclareMathOperator*{\argmin}{argmin}
\usepackage{multicol}

\title{A Majorization-Minimization-Based Method for Nonconvex Inverse Rig Problems in Facial Animation: Algorithm Derivation}

\date{January 2023}	

\author{ \href{https://orcid.org/0000-0002-5656-9189}{\includegraphics[scale=0.06]{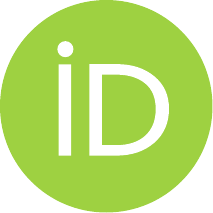}\hspace{1mm}Stevo Racković} \\
	Department of Mathematics\\
	Instituto Superior Técnico\\
	Lisbon, Portugal \\
	\texttt{stevo.rackovic@tecnico.ulisboa.pt} \\
	\And
	\href{https://orcid.org/0000-0003-3071-6627}{\includegraphics[scale=0.06]{orcid.pdf}\hspace{1mm}Cláudia Soares} \\
	Department of Computer Sceince\\
	NOVA School of Science and Technology\\
	Caparica, Portugal \\
	\And
	\href{https://orcid.org/0000-0003-3497-5589}{\includegraphics[scale=0.06]{orcid.pdf}\hspace{1mm}Dušan Jakovetić} \\
	Department of Mathematics\\
	University of Novi Sad\\
	Novi Sad, Serbia \\
 	\And
	Zoranka Desnica \\
	3Lateral Animation Studio\\
	Epic Games Company \\
}

\renewcommand{\shorttitle}{ }

\hypersetup{
pdftitle={A Majorization-Minimization-Based Method for Nonconvex Inverse Rig Problems in Facial Animation: Algorithm Derivation},
pdfauthor={Stevo Racković, Cláudia Soares, Dušan Jakovetić, Zoranka Desnica},
pdfkeywords={Inverse Rig, Quadratic blendshape model, Majorization Minimization},
}

\begin{document}
\maketitle

\begin{abstract}
    Automated methods for facial animation are a necessary tool in the modern industry since the standard blendshape head models consist of hundreds of controllers, and a manual approach is painfully slow. Different solutions have been proposed that produce output in real-time or generalize well for different face topologies. However, all these prior works consider a linear approximation of the blendshape function and hence do not provide a high-enough level of detail for modern realistic human face reconstruction. A second-order blendshape approximation leads to higher fidelity facial animation but generates a non-linear least squares optimization problem with high dimensionality. We derive a method for solving the inverse rig in blendshape animation using quadratic corrective terms, which increases accuracy. At the same time, due to the proposed construction of the objective function, it yields a sparser estimated weight vector compared to the state-of-the-art methods. The former feature means lower demand for subsequent manual corrections of the solution, while the latter indicates that the manual modifications are also easier to include. Our algorithm is iterative and employs a Majorization-Minimization paradigm to cope with the increased complexity produced by adding corrective terms. The surrogate function is easy to solve and allows for further parallelization on the component level within each iteration. This paper is complementary to an accompanying paper \cite{rackovic2022accurate} where we provide detailed experimental results and discussion, including highly-realistic animation data, and show a clear superiority of the results compared to the state-of-the-art methods.
\end{abstract}

\keywords{Inverse Rig \and Quadratic blendshape model \and Majorization Minimization}


\section{Introduction}
\label{sec:intro}

Human face animation is an increasingly popular field of research within the applied mathematics community, of interest not only for the production of movies and video games but also in virtual reality, education, communication, etc. A common approach to face animation is using blendshapes \cite{lewis2014practice, ccetinaslan2016position} --- blendshapes are a vector representation of the face, where each vector $\textbf{b}_1,...,\textbf{b}_m\in\mathbb{R}^{3n}$ stands for a single atomic deformation of the face, and a resting face is represented by $\textbf{b}_0\in\mathbb{R}^{3n}$. More complex expressions are obtained by combining the blendshape vectors, commonly using a linear delta blendshape function:
\begin{equation}\label{eq:intro_model}
f(w_1,...,w_m) = \textbf{b}_0 + \sum_{i=1}^mw_i\Delta\textbf{b}_i,
\end{equation}
where $\Delta\textbf{b}_i=\textbf{b}_i-\textbf{b}_0$ for $i=1,...,m$, and $w_i\in[0,1]$ are activation weights corresponding to each blendshape. These models provide intuitive controls, even though building the base shapes is demanding in terms of time and effort \cite{neumann2013sparse, wang2020facial, li2013realtime, zhang2020facial}. Inverse rig estimation is a common problem that consists of choosing the right activation weights $w_1,...,w_m$ from (\ref{eq:intro_model}) to produce a predefined expression $\widehat{\textbf{b}}\in\mathbb{R}^{3n}$. It is one of the aspects of blendshape animation that can be automated, hence it is often addressed in the literature, usually posed in a least-squares fashion as 
\begin{equation}
\minimize_{w_1,...,w_m} \|f(w_1,...,w_m)- \widehat{\textbf{b}}\|^2,
\end{equation}
with possibly additional constraints or regularization. Here, and throughout the paper, $\|\cdot\|$ stands for the $l_2$ norm. Possible approaches for solving the inverse rig can be classified into data-based and model-based techniques, where the first group demands large amounts of animation data for training purposes, and the second group works with only a blendshape function and the basis vectors. While various machine learning techniques show excellent performance \cite{holden2016learning, bailey2020fast, song2020accurate, Kim2021DeepLU, deng2006animating, Song2011CharacteristicFR, seol2014tuning, holden2015learning, feng2008realtime, yu2014regression, buoaziz2013online}, such methods demand long animated sequences that cover all the regular facial expressions to train a good regressor. This often makes them infeasible or unprofitable. Conversely, model-based approaches \cite{choe2001analysis, sifakis2005automatic, ccetinaslan2016position, li2010example} rely on applying optimization techniques and do not demand training data. Yet, a precise rig function or a good approximation is necessary to provide high-quality results. Without exception, all the papers that propose model-based solutions work with a linear blendshape function, which does not offer high-enough fidelity for realistic animated faces. 

We proposed a new model-based method for solving the inverse rig problem such that it includes the quadratic corrective terms, which leads to higher accuracy of the fit compared to the standard linear rig approximation \cite{holden2016learning, song2017sparse}. Our method utilizes a common framework of Majorization-Minimization \cite{zhang2007surrogate, lange2000optimization}.

\subsection{Contributions}

In the companion paper \cite{rackovic2022accurate}, we present a novel method for solving the inverse rig problem when the blendshape model is assumed to be quadratic. This method targets a specific subdomain of facial animation --- highly-realistic human face models used in movie and video games production. Here the accurate fit plays a more critical role than the execution time. For this reason, the added complexity of a quadratic blendshape rig is justified since it significantly increases the mesh fidelity of the result. Besides increasing the mesh accuracy, our solution yields fewer activated components than the state-of-the-art methods, which is another desirable property in production. The main contributions of the current paper are to provide a detailed derivation of the proposed inverse rig method and describe results on its convergence guarantees. We refer to \cite{rackovic2022accurate} for further practical and implementation aspects, as well as extensive numerical results on real-world animation data.

The rest of the paper is organized as follows. Section \ref{sec:preliminaries} formulates the inverse rig problem and covers the existing solutions from the literature. Section \ref{sec:proposed_solution} introduces our algorithm and gives a detailed derivation of each step. Section \ref{sec:results} discusses numerical results. Finally, we conclude the paper in Section \ref{sec:discussion}.


\section{Problem formulation and preliminaries}\label{sec:preliminaries}


The main components of the blendshape model are the \textit{neutral mesh} $\textbf{b}_0\in \mathbb{R}^{3n}$ sculpted by an artist, as well as a set of $m$ blendshapes $\textbf{b}_1, ..., \textbf{b}_m\in\mathbb{R}^{3n}$, where $n$ is the number of vertices in the mesh. Blendshapes are topologically identical copies of a neutral mesh but with some vertices displaced in space to simulate a local deformation of the face. The offset between neutral mesh and blendshapes yields delta blendshapes $\Delta \textbf{b}_i = \textbf{b}_0 + \textbf{b}_i$, for $ i=1,...,m,$ that are added on top of a neutral mesh $\textbf{b}_0$, with a weight $w_i\in[0,1]$, to produce an effect of local deformation as $ \textbf{b}_0 + w_i\Delta\textbf{b}_i.$ Multiple local deformers are usually combined to produce more complex facial expressions. A blendshape function can then be defined as 
\begin{equation}
f_L(\textbf{w}) = \textbf{b}_0 + \sum_{i=1}^mw_i\Delta\textbf{b}_i = \textbf{b}_0+\textbf{Bw},
\end{equation}
where $\textbf{w} = [w_1,...,w_m]^T$ is a weight vector and $\textbf{B} = [\Delta\textbf{b}_1,...,\Delta\textbf{b}_m]$ is a blendshape matrix. The notation $f_L(\cdot)$ indicates that this is a linear model, while for realistic face representation, it is common to consider a more complex form with quadratic terms, as explained in the following. 

Some pairs of blendshapes, $\textbf{b}_i$ and $\textbf{b}_j$, with an overlapping region of influence might produce artifacts on the face (mesh breaking or giving an unbelievable deformation), hence a corrective term $\textbf{b}^{\{i,j\}}$ needs to be included to fix any issues and make the character appearance natural. An artist discovers these combinations, and corrective terms are conventionally sculpted by hand. Now, whenever the two blendshapes are activated simultaneously, the corrective term is added as well, with a weight equal to the product of two corresponding weights. We define a set $\mathcal{P}$ that stores tuples of indices $(i,j)$ such that a pair of belndshapes $i$ and $j$ have a corresponding corrective term $\textbf{b}^{\{i,j\}}$. A quadratic blendshape function is then defined as:
\begin{equation}\label{eq:quadratic_rig}
f_Q(\textbf{w}) = b_0 + \sum_{i=1}^m w_i\Delta\textbf{b}_i + \sum_{(i,j)\in\mathcal{P}}w_iw_j\textbf{b}^{\{i,j\}}.
\end{equation}

In production, it might be essential to solve the inverse problem. That is, considering there is a given mesh $\widehat{\textbf{b}}$, that is conventionally obtained either as a 3D scan of an actor or a capture of the face markers, one needs to estimate a configuration of the weight vector $\textbf{w}$ that produces a mesh as similar as possible to $\widehat{\textbf{b}}$. This problem is known as the inverse rig problem or the automatic keyframe animation. As we will discuss in the next section, it is common to pose this in a least-squares setting.

\subsection{Existing solutions}\label{sec:existing_sol}

When we consider a model-based solution to the inverse rig problem, the state-of-the-art method is \cite{cetinaslan2020sketching}, where the optimization problem is formulated as regularized least-squares minimization:
\begin{equation}\label{eq:cet}
\minimize_{\textbf{w}}\| \textbf{Bw} -\widehat{\textbf{b}}\|^2 + \alpha\|\textbf{w} \|^2.
\end{equation}
Here, and throughout the paper, $\|\cdot\|$ stands for the $l_2$ norm.
Regularization is necessary because many blendshape pairs are highly correlated, i.e., produce relatively similar deformations over the mesh; hence, the unregularized problem is often ill-posed, and the solution is not unique. Additionally, it is desired to keep the number of non-zero elements of $\textbf{w}$ low because it allows for further manual editing, which is common in animation. The solution to (\ref{eq:cet}) is given in a closed-form as 
\begin{equation}\label{eq:cet_sol}
\textbf{w} = (\textbf{B}^T\textbf{B} + \alpha \textbf{I})^{-1} \textbf{B}^T\widehat{\textbf{b}}.
\end{equation}
A modification to this solution is given in the same paper. Authors approximate a blendshape matrix $\textbf{B}$ with a sparse matrix $\textbf{B}^{loc}$, by applying a heat kernel over the rows of an original blendshape matrix. This sets low values to zero, meaning that effects of the least significant blendshapes are excluded for each vertex. 

A different approach is given in \cite{seol2011artist}, where components are visited and optimized sequentially, and after each iteration $i=1,...,m$, a residual term is updated:
\begin{equation}\label{eq:seol}
\begin{split}
\text{step 1: }\,\,& \minimize_{w_i}  \| w_i\Delta \textbf{b}_i - \widehat{\textbf{b}}\|^2,\\
\text{step 2: }\,\,& \widehat{\textbf{b}}\leftarrow  \widehat{\textbf{b}} - w_i\Delta\textbf{b}_i.
\end{split}
\end{equation}
Here \textit{step 1} finds the optimal activation of a single controller $w_i$, and \textit{step 2} removes its effect for subsequent iterations. This yields a sparse weight vector and excludes the possibility of simultaneously activating mutually exclusive blendshapes (like \texttt{mouth-corner-up} and \texttt{mouth-corner-down}). However, the order in which the blendshapes are visited is crucial to obtain an acceptable solution. The authors propose ordering them based on the average offset that each blendshape causes over the whole face. 

Other papers consider approaches similar to (\ref{eq:cet}) or (\ref{eq:seol}), and they all assume a linear blendshape function. A linear model has an advantage over a quadratic because it gives rise to a convex problem, and it is thus simple and easy to work with; however, that simplicity comes at a price --- it does not provide enough detail for a realistic human face representation in high-quality movies and video games. In the next section, we introduce our solution to the inverse rig problem that takes into account quadratic corrective terms.


\section{Proposed solution}\label{sec:proposed_solution}

This section presents a detailed derivation of our method that utilizes second-order blendshape models; we refer to \cite{rackovic2022accurate} for a detailed method's presentation from the domain point of view and for extensive numerical experiments on the method. Our algorithm targets specifically a high-quality realistic human face animation, hence we assume that a real-time execution is not a priority and that the activation weights are strictly bounded to $[0,1]$ interval\footnote{Many authors neglect this constraint with a justification that the values outside this interval can be still used for \textit{exaggerated} expressions in animated characters. In our setting, this is not allowed by the construction of the models and we have to take these constraints into account.}. We first explain the algorithm derivation and then detail each algorithm step. The optimization problem looks for the optimal weight vector configuration $\textbf{w}$ that fits on a given mesh $\widehat{\textbf{b}}$, assuming a quadratic blendshape function (\ref{eq:quadratic_rig}), as
\begin{equation}\label{eq:qadratic_obj}
\minimize_{\textbf{0}\leq\textbf{w}\leq\textbf{1}} \big\|f_Q(\textbf{w}) - \widehat{\textbf{b}} \big\|^2 + \alpha\textbf{1}^T\textbf{w},
\end{equation}
where the first term is data fidelity and the second is regularization. The non-negativity constraint is important in blendshape animation since negative weights have no semantical meaning and make it harder for animators to adjust the obtained results manually. While weights larger than one might be useful for exaggerated cartoonish expressions, in realistic human avatars this is not a favorable behavior. The regularization term with $\alpha >0$ enforces a low cardinality of the solution vector, which is a desired feature as it makes the results more artist-friendly \cite{seol2011artist}. The problem is approached in a fashion similar to the Levenberg-Marquardt method \cite{ranganathan2004levenberg},  where we choose the initial vector of controller weights $\textbf{w}_{(0)}\in\mathbb{R}^m$, and at each iteration $t=0,...,T$, solve for an optimal increment vector $\textbf{v}\in\mathbb{R}^m$ that solves the following optimization problem:
\begin{equation}\label{eq:inc_qadratic_obj}
	\minimize_{-\textbf{w}_{(t)}\leq\textbf{v}\leq \textbf{1}-\textbf{w}_{(t)}} \big\|f_Q(\textbf{w}_{(t)}+\textbf{v}) - 	\widehat{\textbf{b}} \big\|^2 + \alpha\textbf{1}^T(\textbf{w}_{(t)}+\textbf{v}).
\end{equation}
The weights vector is updated as $\textbf{w}_{(t+1)}=\textbf{w}_{(t)}+\textbf{v}$, and the process is repeated until convergence. Under the quadratic approximation of the rig function, the objective function in (\ref{eq:inc_qadratic_obj}) is fairly complex, hence we simplify it by applying Majorization-Minimization (MM). That is, to solve (\ref{eq:inc_qadratic_obj}), we define an upper bound function $\psi(\textbf{v};\textbf{w}):\mathbb{R}^m\rightarrow \mathbb{R}$ over the objective function in (\ref{eq:inc_qadratic_obj}) such that it is easier to minimize and satisfies \textit{Conditions} 1 and 2 given below. Note that these conditions define a class of functions $\psi(\textbf{v};\textbf{w})$ as potential majorizers to the objective (\ref{eq:inc_qadratic_obj}), but later in this section, we define a specific choice of $\psi(\cdot)$ used in the proposed algorithm.

\textit{Condition 1.} For any feasible vector $\textbf{0}\leq\textbf{w}\leq\textbf{1}$, for all the values of an increment vector $\textbf{v}$ such that $\textbf{0}\leq\textbf{w}+\textbf{v} \leq \textbf{1}$, the following holds:
\begin{equation}
	\|f_Q(\textbf{w+v}) - \widehat{\textbf{b}}\|^2 + \alpha\textbf{1}^T(\textbf{w}+\textbf{v}) \leq \psi(\textbf{v};\textbf{w}).
\end{equation} 

\textit{Condition 2.} The upper bound $\psi(\textbf{v};\textbf{w})$ matches the value of the objective (\ref{eq:inc_qadratic_obj}) at point $\textbf{v}=\textbf{0}$, i.e.,
\begin{equation}
	\|f_Q(\textbf{w}) - \widehat{\textbf{b}}\|^2 + \alpha\textbf{1}^T(\textbf{w}) = \psi(\textbf{0};\textbf{w}).
\end{equation}
In the rest of the paper we will write $\psi(\textbf{v})$ instead of $\psi(\textbf{v};\textbf{w})$, for the sake of simplicity. Further we proceed with the problem in the form 
\begin{equation}\label{eq:increment_problem}
    \minimize_{-\textbf{w}\leq\textbf{v}\leq\textbf{1-w}} \psi(\textbf{v}),
\end{equation}
and the following proposition gives us guarantees that such an approach leads to the minimization of the original objective. \\

\textit{Proposition 1}. Under \textit{Conditions} 1 and 2, a sequence of iterates $\textbf{w}_{(t)}$ for $t\in\mathbb{N}$ obtained as the solutions to problem (\ref{eq:increment_problem}) (with $\textbf{w}_{(t+1)} = \textbf{w}_{(t)}+\textbf{v}_{(t)}$) produces a monotonically non-increasing sequence of values of the objective (\ref{eq:qadratic_obj}), i.e.,
\begin{equation}
	\|f_Q(\textbf{w}_{(t)}) - \widehat{\textbf{b}} \|^2 + \alpha\textbf{1}^T\textbf{w}_{(t)} \geq \|f_Q(\textbf{w}_{(t+1)}) - \widehat{\textbf{b}} \|^2 + \alpha\textbf{1}^T\textbf{w}_{(t+1)} \,\, , \, t\in\mathbb{N}. 
\end{equation}

\begin{proof}
	If $\textbf{v}_{(t)}$ is a minimizer of (\ref{eq:increment_problem}) at iteration $t$, i.e., $$\textbf{v}_{(t)} = \argmin_{-\textbf{w}_{(t)}\leq\textbf{v}\leq\textbf{1-w}_{(t)}} \psi(\textbf{v}),$$
	then we have $\psi(\textbf{0}) \geq \psi(\textbf{v}_{(t)})$. From this and from \textit{Conditions} 2 and 1, we have the following relation:
	\begin{equation}
		\begin{split}
		\|f_Q(\textbf{w}_{(t)}) - \widehat{\textbf{b}}\|^2 & + \alpha\textbf{1}^T(\textbf{w}_{(t)}) \\ 
		= & \psi(\textbf{0}) \geq \psi(\textbf{v}_{(t)}) \\ 
		\geq & \|f_Q(\textbf{w}_{(t)} + \textbf{v}_{(t)}) - \widehat{\textbf{b}}\|^2 + \alpha\textbf{1}^T(\textbf{w}_{(t)}+\textbf{v}_{(t)}),
		\end{split}
	\end{equation}
	which proves the proposition. 
\end{proof}

Additionally, from \cite{wu1983convergence}, we have that under the above iterative method, the objective function converges to a local optimum or a saddle point as number of iterations $t$ goes to infinity.

\paragraph{Upper Bound Function.}

Let us now introduce a specific majorizer that we apply in the proposed algorithm, and below we will give a complete derivation. We define the upper bound function of objective (\ref{eq:inc_qadratic_obj}) as 
\begin{equation}\label{eq:bound_psi}
    \psi(\textbf{v}) = \sum_{i=1}^{3n} \psi_i(\textbf{v}) + \alpha\textbf{1}^T(\textbf{w}+\textbf{v}),
\end{equation}
which is an original regularization term added on a sum of coordinate-wise upper bounds $\psi_i(\textbf{v}):\mathbb{R}^m\rightarrow\mathbb{R}$ of the data fidelity term, where $n$ is the number of vertices in the mesh. The component-wise bounds have the form 
\begin{equation}\label{eq:component_wise_psi}
\psi_i(\textbf{v}) := g_i^2 + 2g_i\sum_{j=1}^m h_{ij}v_j +  
2\big(g_i\lambda_M( \textbf{D}^{(i)},g_i) + \|\textbf{h}_i\|^2\big)\sum_{j=1}^mv_j^2 + 2m\sigma^2(\textbf{D}^{(i)})\sum_{j=1}^mv_j^4
\end{equation}
where $g_i:=\textbf{B}_i\textbf{w} + \textbf{w}^T\textbf{D}^{(i)}\textbf{w} - \widehat{b}_i$, and $\textbf{h}_i := \textbf{B}_i + 2\textbf{w}^T\textbf{D}^{(i)}$ are introduced to simplify the notation; $\textbf{D}^{(i)}\in\mathbb{R}^{m\times m}$ is a symmetric (and sparse) matrix whose nonzero entries are extracted from the corrective blendshapes as $D^{(i)}_{jk} = D^{(i)}_{kj} = \frac{1}{2}b_i^{\{j,k\}}$; the largest singular value of a matrix $\textbf{D}^{(i)}$ is denoted $\sigma(\textbf{D}^{(i)})$; a function $\lambda_M(\textbf{D}^{(i)},g_i):(\mathbb{R}^{m\times m},\mathbb{R})\rightarrow\mathbb{R}$ is defined as
$$  \lambda_M(\textbf{D}^{(i)},g_i) :=
    \begin{cases}
      \lambda_{\text{min}}(\textbf{D}^{(i)}) & \text{if } g_i<0,\\
      \lambda_{\text{max}}(\textbf{D}^{(i)}) & \text{if } g_i\geq0,
    \end{cases} 
$$
where $\lambda_{\text{min}}(\textbf{D}^{(i)})$ represents the smallest and $\lambda_{\text{max}}(\textbf{D}^{(i)})$ the largest eigenvalue of $\textbf{D}^{(i)}$. Under the surrogate function (\ref{eq:bound_psi}), the problem (\ref{eq:increment_problem}) can be analytically solved component-wise, where for each component $j=1,...,m$, the objective is a scalar quartic equation:
\begin{equation}\label{eq:increment_vector}
\begin{split}
    \minimize_{v_j} \,\, qv_j + rv_j^2 + sv_j^4, \\ 
    \text{s.t.   } \, 0\leq w_j+v_j \leq1.
\end{split}
\end{equation}
The expressions for the polynomial coefficients $q,r,s$ are 
\begin{equation}\label{eq:polinomial_coeficients}
q=  2\sum_{i=1}^{3n} g_ih_{ij} + \alpha, \,\, r=  2\sum_{i=1}^{3n}(g_i\lambda_M(\textbf{D}^{(i)},g_i)+\|\textbf{h}_i\|^2),\,\, s=  2m\sum_{i=1}^{3n}\sigma^2(\textbf{D}^{(i)}) .
\end{equation}
Notice that the coefficient $q$ depends on a coordinate $j$, so it has to be computed for each controller separately, while $r$ and $s$ are calculated only once per iteration. We can find the extreme values of the polynomial using the roots of the cubic derivative $ q + 2rv_j + 4sv_j^3 = 0,$ and, if they are within the feasible interval $[0, 1]$, compare them with the polynomial values at the borders to get the constrained minimizer. The idea of component-wise optimization of the weight vector makes our approach somewhat similar to that of \cite{seol2011artist}, yet we do not update the vector until all the components are estimated. This solves the issue of the order in which the controllers are visited and additionally gives a possibility for a parallel implementation of the procedure. Another difference is that \cite{seol2011artist} considers only a single run over the coordinates, while for us, this is only a single step --- we refer to it as an \textit{inner iteration}. In this sense, we start with an arbitrary and feasible vector and repeat the \textit{inner iteration} multiple times to provide an increasingly good estimate of the solution to (\ref{eq:inc_qadratic_obj}), in a manner similar to the trust-region method \cite{tarzanagh2015new, porcelli2013convergence}. A complete inner iteration is summarized in Algorithm \ref{alg:inner_iteration}. In the following lines, we will detail a complete derivation of the upper bound given in (\ref{eq:bound_psi}).

\paragraph{Derivation of the Upper Bound.}

If we write down the data fidelity term from (\ref{eq:inc_qadratic_obj}) as a sum, and consider each element $i=1,...,3n$ of the sum separately, we can represent it in a canonical quadratic form:
\begin{equation}\label{eq:derivation_1}
\sum_{i=1}^{3n}\bigg(\textbf{B}_i\textbf{w} + \sum_{(j,k)\in\mathcal{P}}w_jw_kb_i^{\{j,k\}}-\widehat{{b}}_i\bigg)^2  =  \sum_{i=1}^{3n}\left(\textbf{B}_i\textbf{w}+ \textbf{w}^T\textbf{D}^{(i)}\textbf{w} -\widehat{{b}}_i\right)^2.
\end{equation}
Define a function $\phi_i(\textbf{w}): \mathbb{R}^m\rightarrow\mathbb{R}$ as $
\phi_i(\textbf{w}) : = (\textbf{B}_i\textbf{w}+\textbf{w}^T\textbf{D}^{(i)}\textbf{w}-\widehat{{b}}_i)^2$, which is a single coordinate of the data fidelity term. When we add the increment vector $\textbf{v}$ on top of the current weight vector $\textbf{w}$, this yields:
\begin{equation}\label{eq:phi_i_inc}
\begin{split}
\phi_i(\textbf{w}+\textbf{v}) & = \left(g_i + \textbf{h}_i\textbf{v} + \textbf{v}^T\textbf{D}^{(i)}\textbf{v}\right)^2 \\
& = g_i^2 + 2g_i\textbf{h}_i\textbf{v} + 2g_i\textbf{v}^T\textbf{D}^{(i)}\textbf{v} + \left(\textbf{h}_i\textbf{v} + \textbf{v}^T\textbf{D}^{(i)}\textbf{v}\right)^2.
\end{split}
\end{equation}
The data fidelity term is a sum of functions $\phi_i(\textbf{w})$, hence in order to bound the objective, we bound each element of the sum $\phi_i(\textbf{w}+\textbf{v})\leq\psi_i(\textbf{v})$. Let us first consider two nonlinear terms of $\phi_i(\textbf{w})$, and bound each of them separately. A bound over the quadratic term $2g_i\textbf{v}^T\textbf{D}^{(i)}\textbf{v}$ depends on the sign of $g_i$:
\begin{equation}
2g_i\textbf{v}^T\textbf{D}^{(i)}\textbf{v} \leq 2g_i\lambda_M\big(\textbf{D}^{(i)},g_i\big)\|\textbf{v}\|^2.
\end{equation}
The bound over a quartic term $(\textbf{h}_i\textbf{v} + \textbf{v}^T\textbf{D}^{(i)}\textbf{v})^2$ is obtained by applying the Cauchy-Schwartz inequality three times:
\begin{equation}\label{eq:derivation_4}
\begin{split}
(\textbf{h}_i\textbf{v} + \textbf{v}^T\textbf{D}^{(i)}\textbf{v})^2 & \leq 2(\textbf{h}_i\textbf{v})^2 + 2(\textbf{v}^T\textbf{D}^{(i)}\textbf{v})^2\\
\leq  2\|\textbf{h}_i\|^2\|\textbf{v}\|^2 + 2 \|\textbf{v}\|^4\|\textbf{D}^{(i)}\|^2 & \leq 2\|\textbf{h}_i\|^2\|\textbf{v}\|^2 + 2m\sigma^2(\textbf{D}^{(i)})\sum_{j=1}^mv_j^4.
\end{split}
\end{equation}
The component-wise bound $\psi_i(\textbf{v})$ is then given by (\ref{eq:component_wise_psi}) and a complete upper bound is obtained by summing the component-wise bounds and adding the regularization term (Eq. \ref{eq:bound_psi}).

\begin{algorithm}
\caption{Inner Iteration}\label{alg:inner_iteration}
\begin{algorithmic}
\Require Blendshape matrix $\textbf{B}\in\mathbb{R}^{3n\times m}$, corrective blendshape matrices $\textbf{D}^{(i)}\in\mathbb{R}^{m\times m}$ for $i=1,...,3n$, target mesh $\widehat{\textbf{b}}\in\mathbb{R}^{3n}$, regularization parameter $\alpha\geq0$ and weight vector $\textbf{w}\in[0,1]^m$.
\Ensure $\hat{\textbf{v}}$ - an optimal increment vector as a solution to (\ref{eq:increment_vector}).
\State Compute coefficients $q,r$,$s$ from Eq. (\ref{eq:polinomial_coeficients}) and solve for an optimal increment vector $\hat{\textbf{v}}$:
\State $r=  2\sum_{i=1}^{3n}(g_i\lambda_M(\textbf{D}^{(i)},g_i)+\|\textbf{h}_i\|^2)$, 
\State $s=  2m\sum_{i=1}^{3n}\sigma^2(\textbf{D}^{(i)})$, 
\For{$k=1,...,m$} 
    \State $q=  2\sum_{i=1}^{3n} g_ih_{ij} + \alpha$
    \State $\hat{v}_k = \argmin_{v} qv + rv^2 + sv^4 $
    \State \qquad s.t. $-w_k\leq v \leq 1-w_k$
\EndFor
\State \Return $\hat{\textbf{v}}$
\end{algorithmic}
\end{algorithm}

\paragraph{Feasibility.}

Let us now state another proposition showing that the above derived upper bound function is feasible for \textit{Proposition 1}.\\

\textit{Proposition 2}. The surrogate function $\psi(\textbf{v}):\mathbb{R}^m\rightarrow \mathbb{R}$, defined as in (\ref{eq:bound_psi}) satisfies \textit{Conditions 1} and \textit{2}, and hence it satisfies \textit{Proposition 1}.
\begin{proof}
In (\ref{eq:phi_i_inc})-(\ref{eq:derivation_4}), function $\psi(\cdot)$ is derived so that it bounds the objective function (\ref{eq:derivation_1}) from above, hence it satisfies \textit{Condition 1} by construction.

Now recall that the data fidelity term is a sum of functions $\phi_i(\cdot)$. To prove that \textit{Condition 2} is satisfied, it suffices to show that $\phi_i(\textbf{w}+\textbf{0})=\psi_i(\textbf{0})$. From (\ref{eq:phi_i_inc}) we see that $\phi_i(\textbf{w}+\textbf{0}) = g_i^2$, and from (\ref{eq:component_wise_psi}) we get $\psi_i(\textbf{0}) = g_i^2$, which proves it. Hence, the derived bound satisfies \textit{Proposition 1}. 
\end{proof}

\paragraph{Complete Algorithm.}

As mentioned above, our solution is iterative, and based on applying Algorithm \ref{alg:inner_iteration} until convergence. In each iteration $t$, the weight vector is updated by adding an estimated increment vector: $\textbf{w}_{(t+1)} = \textbf{w}_{(t)}+\textbf{v}.$ The algorithm terminates when either a given maximum number of iterations is reached, or the difference between the cost function in two consecutive iterations is below a given threshold value $\epsilon>0$. While any feasible vector $\textbf{0}\leq\textbf{w}\leq\textbf{1}$ can be used for initialization, we can rely on domain knowledge to choose a good starting point leading to faster convergence, and these strategies are discussed in the companion paper \cite{rackovic2022accurate}. A complete method is summarized in Algorithm \ref{alg:the_algorithm}. Notice that in one of the steps of the algorithm, we compute eigen- and singular values of matrices $\textbf{D}^{(i)}$. This computation is needed only once per animated character, and we can reuse the calculated values for each following frame that is to be fitted.

\begin{algorithm}
\caption{Proposed Method}\label{alg:the_algorithm}
\begin{algorithmic}
\Require 
Blendshape matrix $\textbf{B}\in\mathbb{R}^{3n\times m}$, corrective blendshapes $\textbf{b}^{\{i,j\}}\in\mathbb{R}^{3n}$ for $(i,j)\in\mathcal{P}$, target mesh $\widehat{\textbf{b}}\in\mathbb{R}^{3n}$, regularization parameter $\alpha\geq0$, initial weight vector $\textbf{w}_{(0)}\in[0,1]^m$, maximum number of iterations $T\in\mathbb{N}$, tolerance $\epsilon>0$.
\Ensure $\hat{\textbf{w}}$ - an approximate minimizer of the problem (\ref{eq:qadratic_obj}).
\State For each coordinate $i=1,...,3n$ compose a matrix $\textbf{D}^{(i)}\in\mathbb{R}^{m\times m}$ from the corrective terms, and extract singular and eigen values ($\sigma,\lambda_{\text{min}}, \lambda_{\text{max}}$):
\For {$i = 1,...,3n$}
    \For {$(j,k)\in\mathcal{P}$}
        \State $D^{(i)}_{jk} = D^{(i)}_{kj} = 1/2 b^{\{j,k\}}_i.$
    \EndFor
    \State  $\textbf{D}^{(i)} \rightarrow \lambda_{\text{min}}(\textbf{D}^{(i)})$, $\lambda_{\text{max}}(\textbf{D}^{(i)})$, $\sigma(\textbf{D}^{(i)})$.
\EndFor
\State Repeat Algorithm 1 until convergence:
\For{$t = 0,...,T$}
    \State Compute optimal increment $\hat{\textbf{v}}$ using Algorithm \ref{alg:inner_iteration} 
    \State Update the weight vector $\textbf{w}_{(t)}$:
    \State $\textbf{w}_{(t+1)} = \textbf{w}_{(t)} + \hat{\textbf{v}}$
    \State Check convergence
    \If{$\left| \psi(\hat{\textbf{v}}) - \psi(\textbf{0}) \right|<\epsilon$} 
        \State $\hat{\textbf{w}} \leftarrow \textbf{w}_{(t+1)}$
        \State \Return $\hat{\textbf{w}}$
    \EndIf
\EndFor
\State $\hat{\textbf{w}} \leftarrow \textbf{w}_{(t+1)}$
\State \Return $\hat{\textbf{w}}$
\end{algorithmic}
\end{algorithm}

\textit{Corollary 1}.
The estimate sequence $\textbf{w}_{(t)}$ produced by Algorithm \ref{alg:the_algorithm} is feasible for problem (\ref{eq:qadratic_obj}) at all iterations $t\in\mathbb{N}$, as long as the initial weight vector $\textbf{w}_{(0)}$ is feasible.

\section{Results}\label{sec:results}

The experiments in this section are performed over an animated avatar \textit{Omar}, available at Metahuman Creator\footnote{unrealengine.com/en-US/metahuman}. The character consists of $m=130$ base blendshapes and $n=4000$ vertices in the face. An optimized choice of regularization parameter $\alpha=5$ is estimated experimentally from the training data, to give a good trade-off between the high accuracy of the mesh fit and a low cardinality of the weight vector. In our experiments, the weight vector is initialized by (\ref{eq:cet_sol}), where the weights outside of the feasible set are clipped to 0 or 1, and then further optimized using our algorithm. While any feasible vector can be used for initialization, this choice has shown quick convergence and precise mesh fit.

\begin{figure}
	\includegraphics[height=0.21\linewidth]{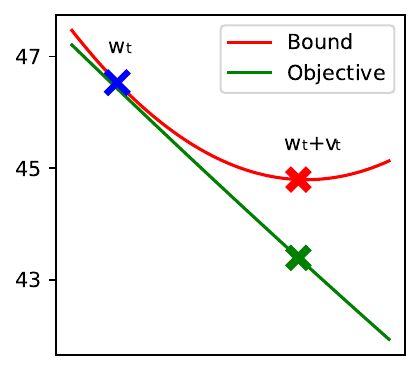}	\includegraphics[height=0.21\linewidth]{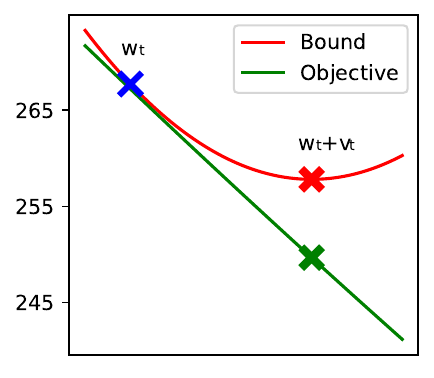}	\includegraphics[height=0.21\linewidth]{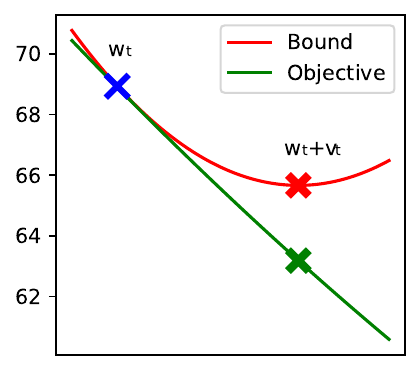}	\includegraphics[height=0.21\linewidth]{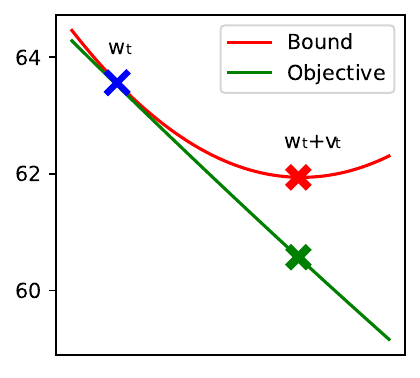}
	\caption{Several example frames showing a relationship between the objective function (green) and the proposed surrogate function (red), at iteration $t$. A current iterate $\textbf{w}_t$ is marked by a blue cross, and the estimated next iterate $\textbf{w}_t + \textbf{v}_t$ is marked by a red cross showing the values of a surrogate function, and by a green one, indicating the actual value of the objective function. }\label{fig:surrogate}
\end{figure}

Let us first examine the relationship between the objective function  (\ref{eq:qadratic_obj}) and its corresponding surrogate function (\ref{eq:bound_psi}). Figure 1 shows several example frames and the behavior of two functions in the initial iterations. The blue cross represents the value of the function at an initial iterate $\textbf{w}_t$, i.e., for $\textbf{v}=0$ (As stated in \textit{Condition 2}, the two functions have the same value in this point). Values of the two functions are then denoted by the corresponding color cross in the estimated iterate $\textbf{w}_t+\textbf{v}_t$. The other values of the two curves are obtained by interpolating the two vectors $\textbf{w}_t$ and $\textbf{w}_t+\textbf{v}_t$ and estimating the corresponding function value. These results clearly show how minimizing the upper bound leads to a nice decrease in the objective. 

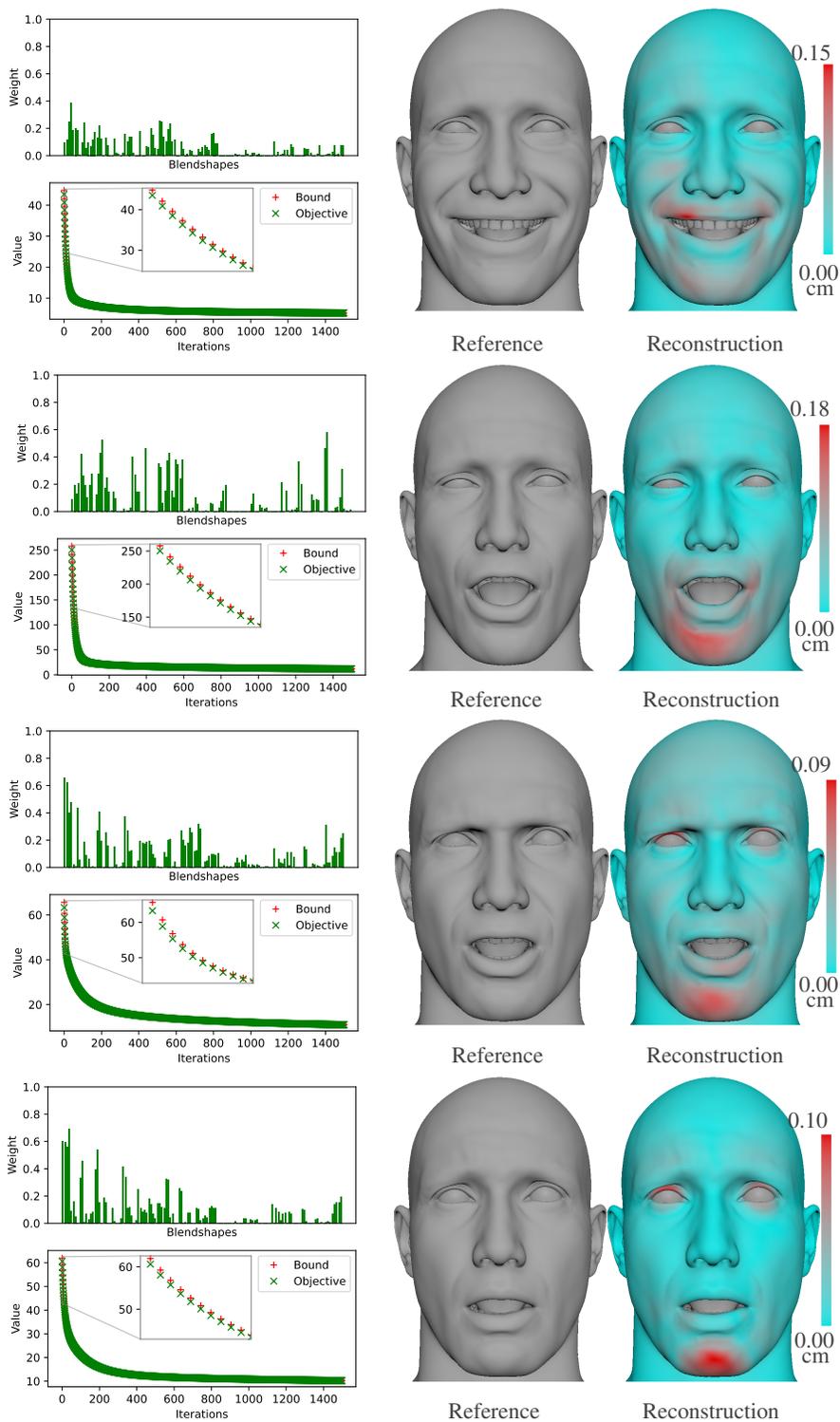
\begin{figure}
	\centering
	\begin{tikzpicture}
		\node[above right, inner sep=0] (image) at (0,0){\includegraphics[height=0.3\textwidth]{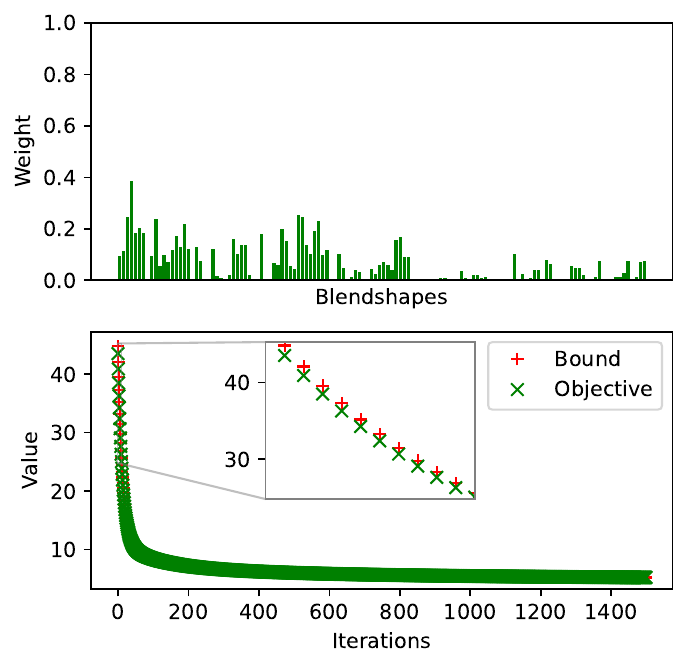}};
		\node[above right, inner sep=0] (image) at (5.5,0){\includegraphics[height=0.3\textwidth]{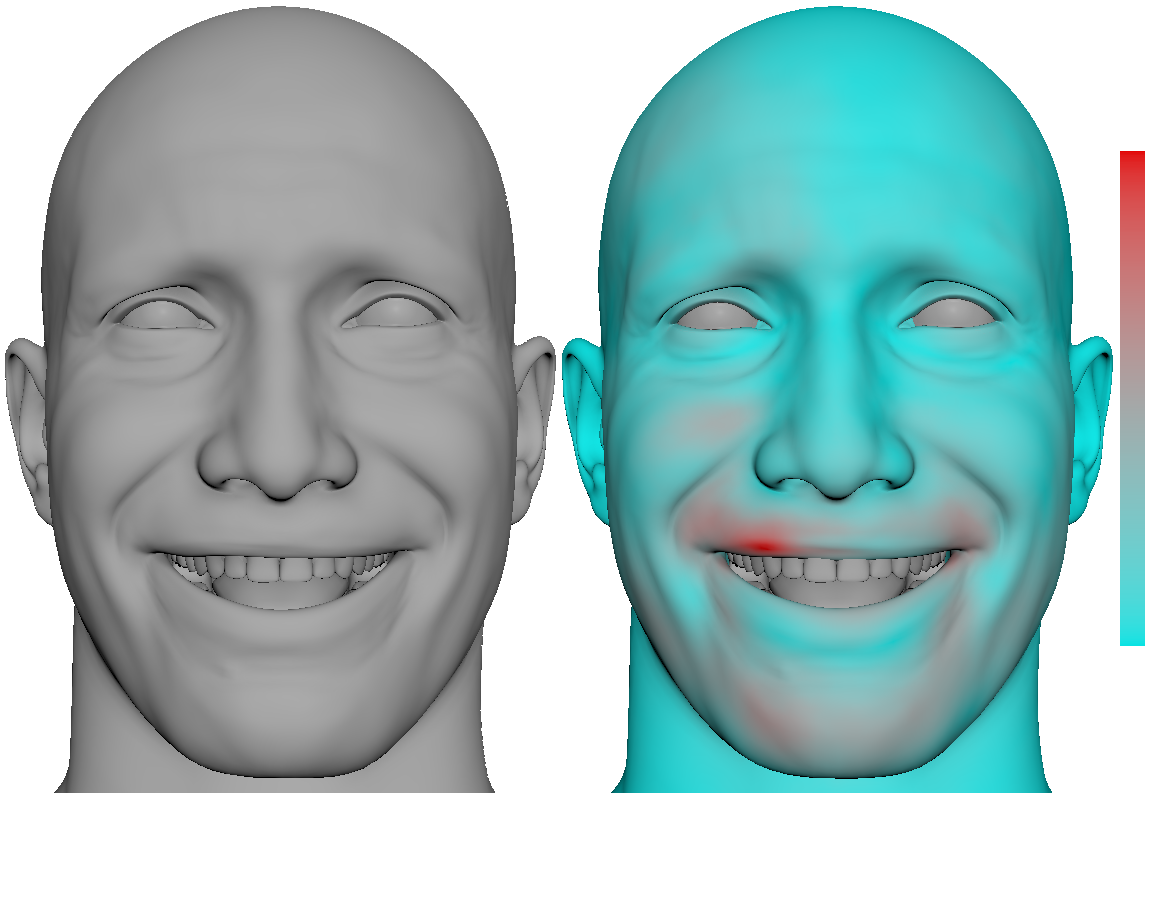}};
		\begin{scope}[
			x={($0.1*(image.south east)$)},
			y={($0.1*(image.north west)$)}]
			\node[darkgray] at (5.70,0.5) {\footnotesize Reference };			
			\node[darkgray] at (8.30,0.5) {\footnotesize Reconstruction };
			\node[darkgray] at (5.6,8.7) {\footnotesize 0.15 };
			\node[darkgray] at (8.6,2.5) {\footnotesize 0.00 };
			\node[darkgray] at (8.85,2.0) {\footnotesize cm };
		\end{scope}
	\end{tikzpicture}
	\begin{tikzpicture}
		\node[above right, inner sep=0] (image) at (0,0){\includegraphics[height=0.3\textwidth]{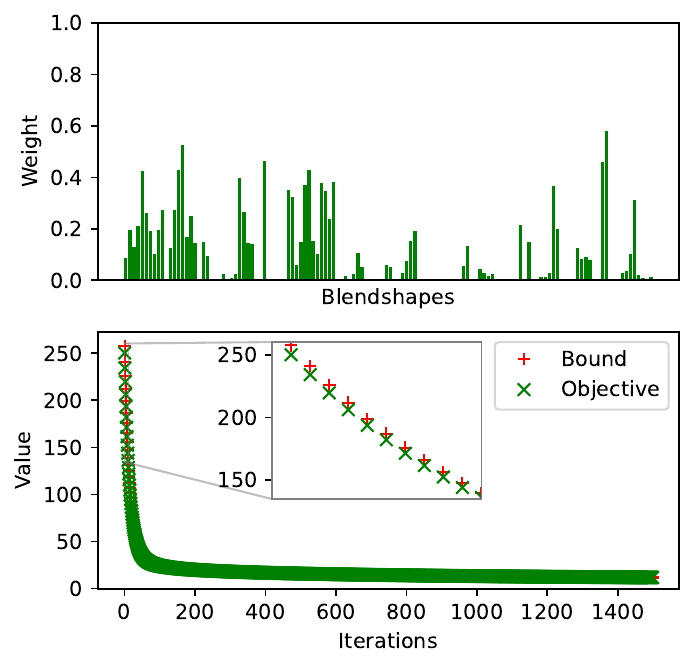}};
		\node[above right, inner sep=0] (image) at (5.5,0){\includegraphics[height=0.3\textwidth]{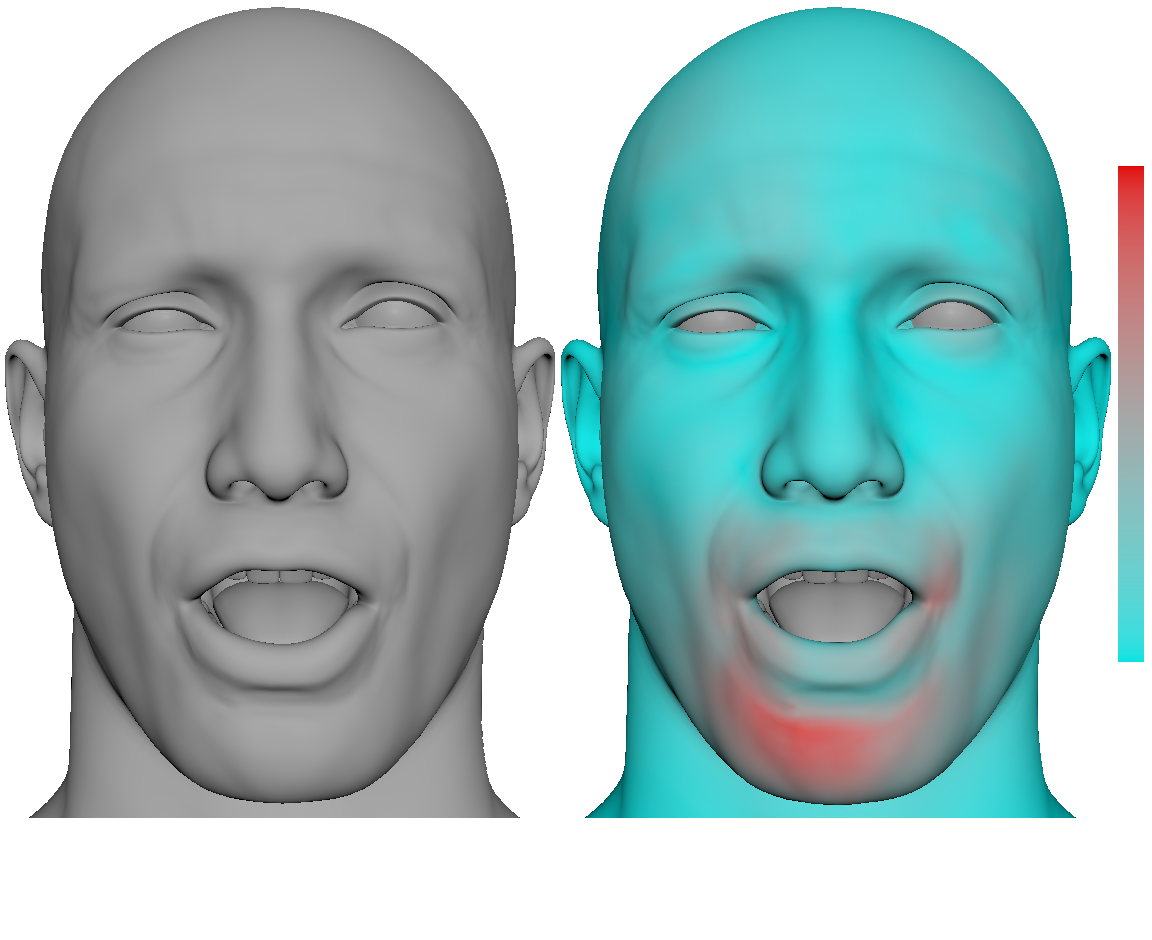}};
		\begin{scope}[
			x={($0.1*(image.south east)$)},
			y={($0.1*(image.north west)$)}]
			\node[darkgray] at (5.70,0.5) {\footnotesize Reference };			
			\node[darkgray] at (8.30,0.5) {\footnotesize Reconstruction };
			\node[darkgray] at (5.6,8.7) {\footnotesize 0.18 };
			\node[darkgray] at (8.6,2.5) {\footnotesize 0.00 };
			\node[darkgray] at (8.85,2.0) {\footnotesize cm };
		\end{scope}
	\end{tikzpicture}
	\begin{tikzpicture}
		\node[above right, inner sep=0] (image) at (0,0){\includegraphics[height=0.3\textwidth]{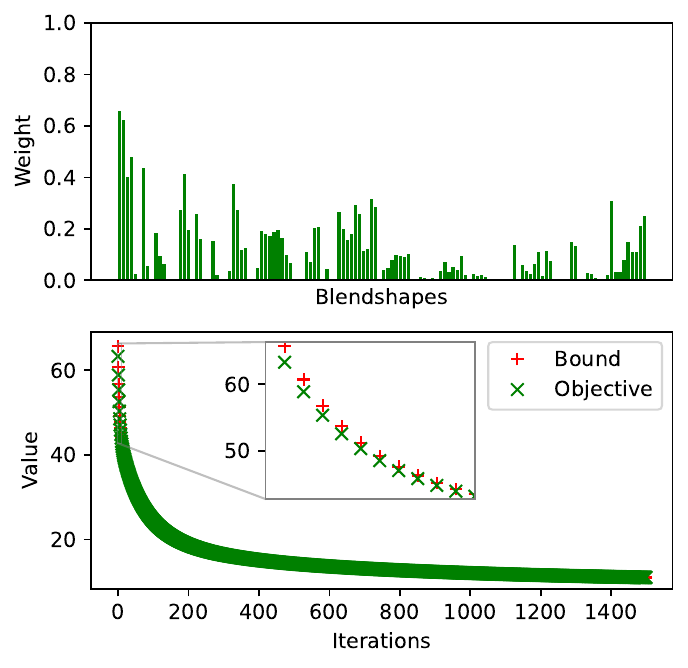}};
		\node[above right, inner sep=0] (image) at (5.5,0){\includegraphics[height=0.3\textwidth]{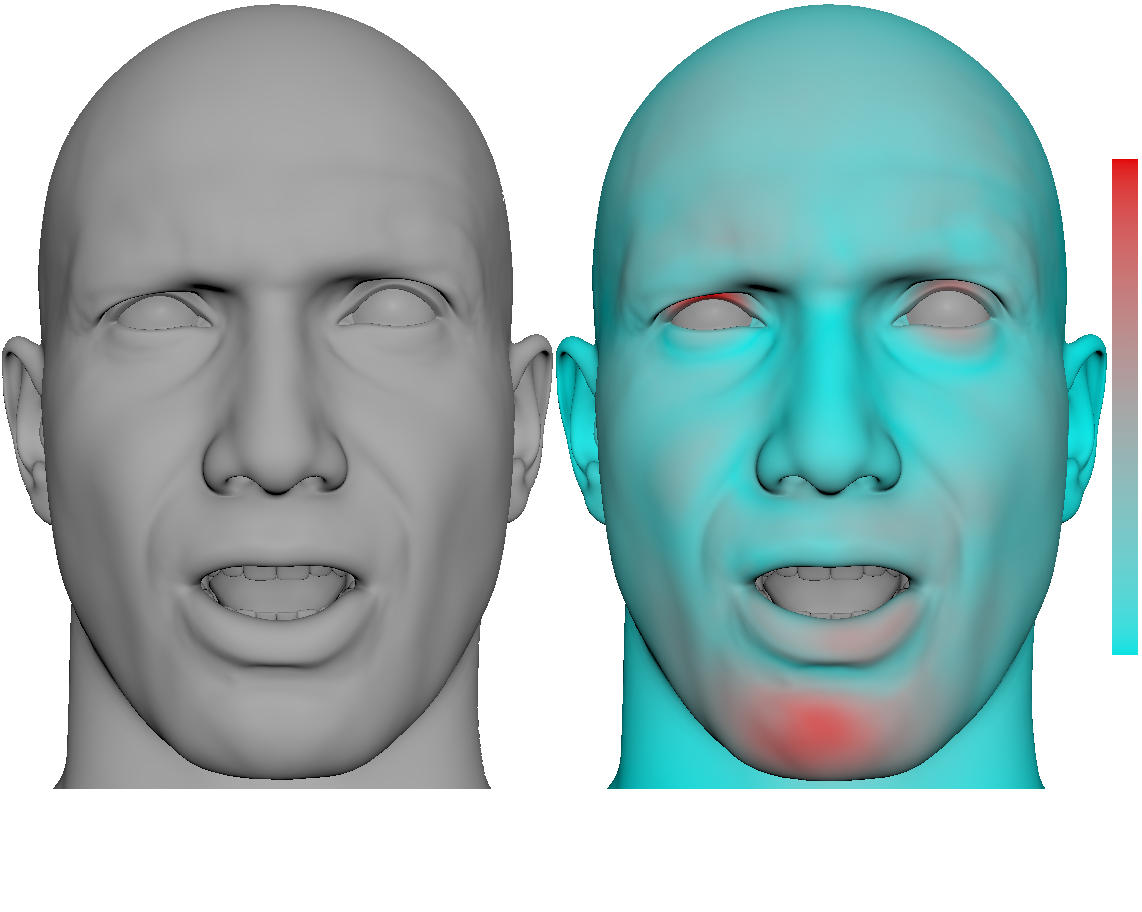}};
		\begin{scope}[
			x={($0.1*(image.south east)$)},
			y={($0.1*(image.north west)$)}]
			\node[darkgray] at (5.70,0.5) {\footnotesize Reference };			
			\node[darkgray] at (8.30,0.5) {\footnotesize Reconstruction };
			\node[darkgray] at (5.6,8.7) {\footnotesize 0.09 };
			\node[darkgray] at (8.6,2.5) {\footnotesize 0.00 };
			\node[darkgray] at (8.85,2.0) {\footnotesize cm };
		\end{scope}
	\end{tikzpicture}
	\begin{tikzpicture}
		\node[above right, inner sep=0] (image) at (0,0){\includegraphics[height=0.3\textwidth]{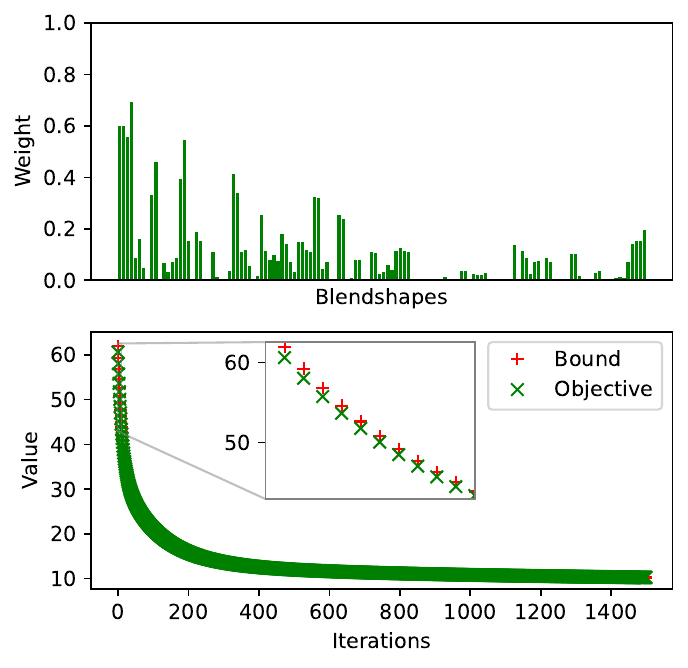}};
		\node[above right, inner sep=0] (image) at (5.5,0){\includegraphics[height=0.3\textwidth]{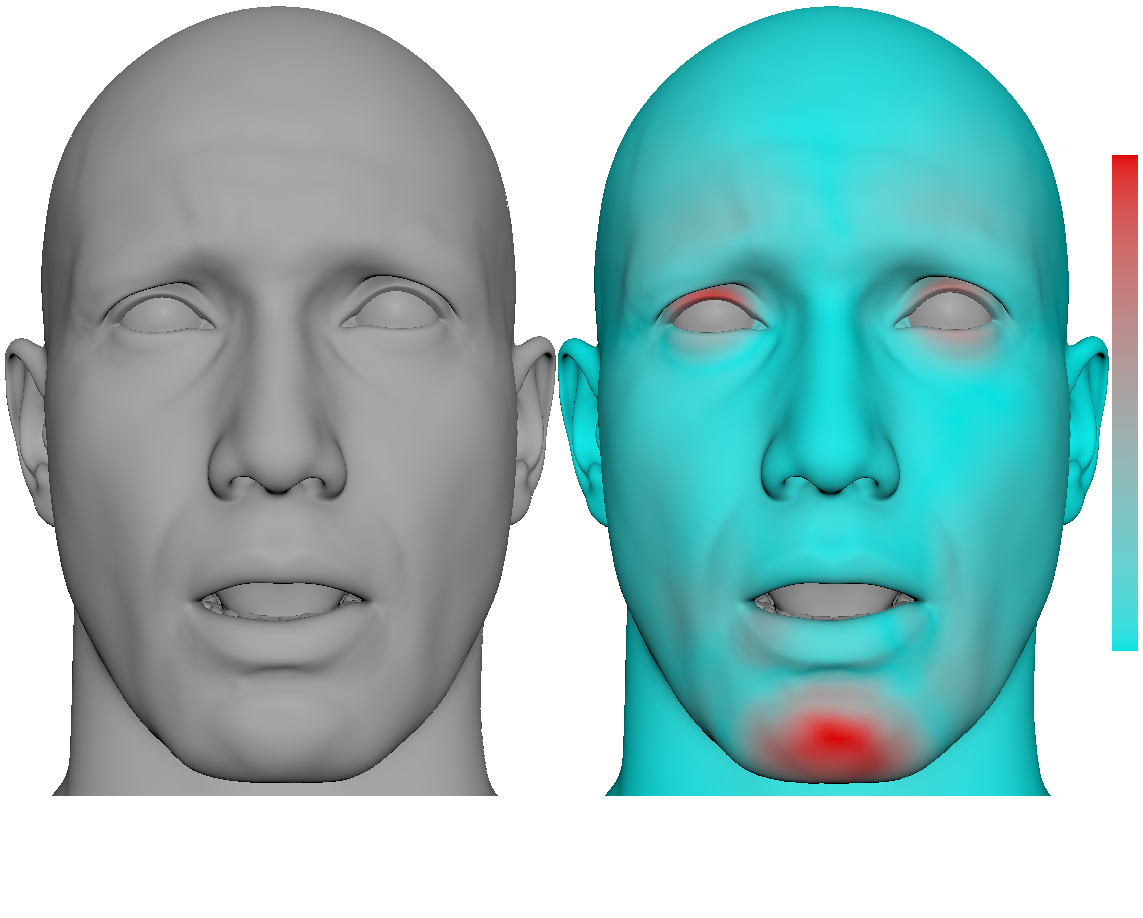}};
		\begin{scope}[
			x={($0.1*(image.south east)$)},
			y={($0.1*(image.north west)$)}]
			\node[darkgray] at (5.70,0.5) {\footnotesize Reference };			
			\node[darkgray] at (8.30,0.5) {\footnotesize Reconstruction };
			\node[darkgray] at (5.6,8.7) {\footnotesize 0.10 };
			\node[darkgray] at (8.6,2.5) {\footnotesize 0.00 };
			\node[darkgray] at (8.85,2.0) {\footnotesize cm };
		\end{scope}
	\end{tikzpicture}
	\caption{Results over four animated frames. Barplots on the left show the estimated activation weight vector $\textbf{w}$; below that is a scatter-plot showing the behavior of the objective function and the upper bound over the iterations; on the right, we have a reference mesh (gray) and the reconstruction by the proposed algorithm --- colors indicate the level of offset over the mesh, corresponding to the color-bar on the right.}\label{fig:res}
\end{figure}

Further, we show the final results of the method in Figure 2. On the left side of the subfigures, one can see a barplot representing the estimated activation weights in the range [0,1] --- we can see that the vector is relatively sparse, and most of the weights have low values. The average cardinality of the weight vector, over 500 test frames, is 85. Below we have the behavior of the objective function over iterations of the algorithm. Green crosses represent the value of the function (\ref{eq:qadratic_obj}) for the estimated iterate $\textbf{w}_t+\textbf{v}_t$, while red pluses denote the value of the corresponding upper bound function (\ref{eq:bound_psi}), in line with the results presented in Figure 1. The plot zooms on initial iterations to show how the proposed bound is tight, and the gap between the surrogate and objective quickly decreases. We can also notice a nice decrease in the objective, confirming the above-stated monotonicity and convergence properties. 

Finally, on the right side of the subfigures, we show the reference mesh $\widehat{\textbf{b}}$ in gray, versus the reconstruction obtained by the proposed method. The reconstructed mesh is colored so that red tones indicate a higher error, according to the color bar on the right. We can see that in each case, the reconstruction is really close to the original frame, and even the regions indicated by the red color are not visually different. The average root mean squared error over 500 test frames is $0.09$.

For a detailed numerical discussion and comparison with other methods see \cite{rackovic2022accurate} --- the reference examines results over multiple animated characters and shows that the proposed method produces accurate and sparse solutions to the inverse rig problem.


\section{Conclusion}
\label{sec:discussion}

This paper introduced the first model-based method for solving the inverse rig problem under the quadratic blendshape function. The proposed algorithm is based on the applications of Majorization-Minimization. A specific surrogate function is derived, and we provide guarantees that it leads to a non-increasing cost sequence of the original non-convex optimization problem (\ref{eq:qadratic_obj}).

The algorithm targets a high-quality facial animation for the video games and movie industry, and hence it assumes that the higher mesh fidelity is more important than the computation speed. For this reason, we rely on the quadratic blendshape function that is more precise than the standard linear one, and also that the weights are strictly constrained to a $[0,1]$ interval to avoid exaggerated or non-credible expressions.

While this paper had the purpose of giving a complete derivation of the proposed algorithm, \cite{rackovic2022accurate} presents an in-depth discussion of the applications and results.


\textbf{Funding}

This work has received funding from the European Union's Horizon 2020 research and innovation program under the Marie Skłodowska-Curie grant agreement No. 812912, from FCT IP strategic project NOVA LINCS (FCT UIDB/04516/2020) and project DSAIPA/AI/0087/2018. The work has also been supported in part by the Ministry of Education, Science and Technological Development of the Republic of Serbia (Grant No. 451-03-9/2021-14/200125).

\bibliographystyle{unsrtnat}
\bibliography{references}  

\begin{thebibliography}{30}
\providecommand{\natexlab}[1]{#1}
\providecommand{\url}[1]{\texttt{#1}}
\expandafter\ifx\csname urlstyle\endcsname\relax
  \providecommand{\doi}[1]{doi: #1}\else
  \providecommand{\doi}{doi: \begingroup \urlstyle{rm}\Url}\fi

\bibitem[Rackovi{\'c} et~al.(2023)Rackovi{\'c}, Soares, Jakoveti{\'c}, and
  Desnica]{rackovic2022accurate}
Stevo Rackovi{\'c}, Cl{\'a}udia Soares, Du{\v{s}}an Jakoveti{\'c}, and Zoranka
  Desnica.
\newblock Accurate and interpretable solution of the inverse rig for realistic
  blendshape models with quadratic corrective terms.
\newblock \emph{arXiv preprint arXiv:2302.04843}, 2023.

\bibitem[Lewis et~al.(2014)Lewis, Anjyo, Rhee, Zhang, Pighin, and
  Deng]{lewis2014practice}
John~P. Lewis, Ken Anjyo, Taehyun Rhee, Mengjie Zhang, Frederic~H. Pighin, and
  Zhigang Deng.
\newblock Practice and theory of blendshape facial models.
\newblock \emph{Eurographics (State of the Art Reports)}, 1\penalty0
  (8):\penalty0 2, 2014.
\newblock \doi{https://doi.org/10.2312/egst.20141042}.

\bibitem[{\c{C}}etinaslan(2016)]{ccetinaslan2016position}
Cumhur~Ozan {\c{C}}etinaslan.
\newblock \emph{Position Manipulation Techniques for Facial Animation}.
\newblock PhD thesis, Faculdade de Ciencias da Universidade do Porto, 2016.

\bibitem[Neumann et~al.(2013)Neumann, Varanasi, Wenger, Wacker, Magnor, and
  Theobalt]{neumann2013sparse}
Thomas Neumann, Kiran Varanasi, Stephan Wenger, Markus Wacker, Marcus Magnor,
  and Christian Theobalt.
\newblock Sparse localized deformation components.
\newblock \emph{ACM Trans. Graph.}, 32\penalty0 (6):\penalty0 1--10, 2013.
\newblock \doi{https://doi.org/10.1145/2508363.2508417}.

\bibitem[Mengjiao et~al.(2020)Mengjiao, Derek, Stefanos, and
  Thabo]{wang2020facial}
Wang Mengjiao, Bradley Derek, Zafeiriou Stefanos, and Beeler Thabo.
\newblock Facial expression synthesis using a global-local multilinear
  framework.
\newblock In \emph{Computer Graphics Forum}, volume~39, pages 235--245. Wiley
  Online Library, 2020.
\newblock \doi{https://doi.org/10.1111/cgf.13926}.

\bibitem[Li et~al.(2013)Li, Yu, Ye, and Bregler]{li2013realtime}
Hao Li, Jihun Yu, Yuting Ye, and Chris Bregler.
\newblock Realtime facial animation with on-the-fly correctives.
\newblock \emph{ACM Trans. Graph.}, 32\penalty0 (4):\penalty0 42--1, 2013.
\newblock \doi{https://doi.org/10.1145/2461912.2462019}.

\bibitem[Zhang et~al.(2020)Zhang, Chen, and Zheng]{zhang2020facial}
Juyong Zhang, Keyu Chen, and Jianmin Zheng.
\newblock Facial expression retargeting from human to avatar made easy.
\newblock \emph{IEEE Transactions on Visualization and Computer Graphics},
  28\penalty0 (2):\penalty0 1274--1287, 2020.
\newblock \doi{https://doi.org/10.1109/TVCG.2020.3013876}.

\bibitem[Holden et~al.(2016)Holden, Saito, and Komura]{holden2016learning}
Daniel Holden, Jun Saito, and Taku Komura.
\newblock Learning inverse rig mappings by nonlinear regression.
\newblock \emph{IEEE Transactions on Visualization and Computer Graphics},
  23\penalty0 (3):\penalty0 1167--1178, 2016.
\newblock \doi{https://doi.org/10.1109/TVCG.2016.2628036}.

\bibitem[Bailey et~al.(2020)Bailey, Omens, Dilorenzo, and
  O'Brien]{bailey2020fast}
Stephen~W. Bailey, Dalton Omens, Paul Dilorenzo, and James~F. O'Brien.
\newblock Fast and deep facial deformations.
\newblock \emph{ACM Trans. Graph.}, 39\penalty0 (4), 2020.
\newblock \doi{https://doi.org/10.1145/3386569.3392397}.

\bibitem[Song et~al.(2020)Song, Shi, and Reed]{song2020accurate}
Steven~L. Song, Weiqi Shi, and Michael Reed.
\newblock Accurate face rig approximation with deep differential subspace
  reconstruction.
\newblock \emph{ACM Trans. Graph.}, 39, 2020.
\newblock \doi{https://doi.org/10.1145/3386569.3392491}.

\bibitem[Seonghyeon et~al.(2021)Seonghyeon, Sunjin, Kwanggyoon, Roger, and
  Junyong]{Kim2021DeepLU}
Kim Seonghyeon, Jung Sunjin, Seo Kwanggyoon, Blanco I.~Ribera Roger, and Noh
  Junyong.
\newblock Deep learning‐based unsupervised human facial retargeting.
\newblock \emph{Computer Graphics Forum}, 40, 2021.
\newblock \doi{https://doi.org/10.1111/cgf.14400}.

\bibitem[Deng et~al.(2006)Deng, Chiang, Fox, and Neumann]{deng2006animating}
Zhigang Deng, Pei-Ying Chiang, Pamela Fox, and Ulrich Neumann.
\newblock Animating blendshape faces by cross-mapping motion capture data.
\newblock In \emph{Proceedings of the 2006 Symposium on Interactive 3D Graphics
  and Games}, page 43–48. Association for Computing Machinery, 2006.
\newblock \doi{https://doi.org/10.1145/1111411.1111419}.

\bibitem[Song et~al.(2011)Song, Choi, Seol, and yong
  Noh]{Song2011CharacteristicFR}
Jaewon Song, Byungkuk Choi, Yeongho Seol, and Jun yong Noh.
\newblock Characteristic facial retargeting.
\newblock \emph{Computer Animation and Virtual Worlds}, 22, 2011.
\newblock \doi{https://doi.org/10.1002/cav.414}.

\bibitem[Seol and Lewis(2014)]{seol2014tuning}
Yeongho Seol and J.~P. Lewis.
\newblock Tuning facial animation in a mocap pipeline.
\newblock In \emph{ACM SIGGRAPH 2014 Talks}, SIGGRAPH '14. Association for
  Computing Machinery, 2014.
\newblock \doi{https://doi.org/10.1145/2614106.2614108}.

\bibitem[Holden et~al.(2015)Holden, Saito, and Komura]{holden2015learning}
Daniel Holden, Jun Saito, and Taku Komura.
\newblock Learning an inverse rig mapping for character animation.
\newblock In \emph{Proceedings of the 14th ACM SIGGRAPH/Eurographics Symposium
  on Computer Animation}, pages 165--173, 2015.
\newblock \doi{https://doi.org/10.1145/2786784.2786788}.

\bibitem[Feng et~al.(2008)Feng, Kim, and Yu]{feng2008realtime}
Wei-Wen Feng, Byung-Uck Kim, and Yizhou Yu.
\newblock Real-time data driven deformation using kernel canonical correlation
  analysis.
\newblock \emph{ACM Trans. Graph.}, 27\penalty0 (3):\penalty0 1–9, 2008.
\newblock \doi{https://doi.org/10.1145/1360612.1360690}.

\bibitem[Yu and Liu(2014)]{yu2014regression}
Hui Yu and Honghai Liu.
\newblock Regression-based facial expression optimization.
\newblock \emph{IEEE Transactions on Human-Machine Systems}, 44\penalty0
  (3):\penalty0 386--394, 2014.
\newblock \doi{https://doi.org/10.1109/THMS.2014.2313912}.

\bibitem[Bouaziz et~al.(2013)Bouaziz, Wang, and Pauly]{buoaziz2013online}
Sofien Bouaziz, Yangang Wang, and Mark Pauly.
\newblock Online modeling for realtime facial animation.
\newblock \emph{ACM Trans. Graph.}, 32\penalty0 (4), 2013.
\newblock \doi{https://doi.org/10.1145/2461912.2461976}.

\bibitem[Choe and Ko(2006)]{choe2001analysis}
Byoungwon Choe and Hyeong-Seok Ko.
\newblock Analysis and synthesis of facial expressions with hand-generated
  muscle actuation basis.
\newblock In \emph{ACM SIGGRAPH 2006 Courses}. 2006.
\newblock \doi{https://doi.org/10.1145/1198555.1198595}.

\bibitem[Sifakis et~al.(2005)Sifakis, Neverov, and
  Fedkiw]{sifakis2005automatic}
Eftychios Sifakis, Igor Neverov, and Ronald Fedkiw.
\newblock Automatic determination of facial muscle activations from sparse
  motion capture marker data.
\newblock In \emph{ACM SIGGRAPH 2005 Papers}, page 417–425. Association for
  Computing Machinery, 2005.
\newblock \doi{https://doi.org/10.1145/1186822.1073208}.

\bibitem[Li et~al.(2010)Li, Weise, and Pauly]{li2010example}
Hao Li, Thibaut Weise, and Mark Pauly.
\newblock Example-based facial rigging.
\newblock \emph{ACM Trans. Graph.}, 29\penalty0 (4):\penalty0 1--6, 2010.
\newblock \doi{https://doi.org/10.1145/1778765.1778769}.

\bibitem[Song et~al.(2017)Song, Blanco~i Ribera, Cho, You, Lewis, Choi, and
  Noh]{song2017sparse}
Jaewon Song, Roger Blanco~i Ribera, Kyungmin Cho, Mi~You, John~P Lewis,
  Byungkuk Choi, and Junyong Noh.
\newblock Sparse rig parameter optimization for character animation.
\newblock In \emph{Computer Graphics Forum}, volume~36, pages 85--94. Wiley
  Online Library, 2017.

\bibitem[Zhang et~al.(2007)Zhang, Kwok, and Yeung]{zhang2007surrogate}
Zhihua Zhang, James Tin-Yau Kwok, and Dit-Yan Yeung.
\newblock Surrogate maximization/minimization algorithms and extensions.
\newblock 2007.

\bibitem[Lange et~al.(2000)Lange, Hunter, and Yang]{lange2000optimization}
Kenneth Lange, David~R Hunter, and Ilsoon Yang.
\newblock Optimization transfer using surrogate objective functions.
\newblock \emph{Journal of computational and graphical statistics}, 9\penalty0
  (1):\penalty0 1--20, 2000.

\bibitem[Cetinaslan and Orvalho(2020)]{cetinaslan2020sketching}
Ozan Cetinaslan and Veronica Orvalho.
\newblock Sketching manipulators for localized blendshape editing.
\newblock \emph{Graphical Models}, 108:\penalty0 101059, 2020.
\newblock \doi{https://doi.org/10.1016/j.gmod.2020.101059}.

\bibitem[Seol et~al.(2011)Seol, Seo, Kim, Lewis, and Noh]{seol2011artist}
Yeongho Seol, Jaewoo Seo, Paul~Hyunjin Kim, J.~P. Lewis, and Junyong Noh.
\newblock Artist friendly facial animation retargeting.
\newblock In \emph{Proceedings of the 2011 SIGGRAPH Asia Conference}.
  Association for Computing Machinery, 2011.
\newblock \doi{https://doi.org/10.1145/2024156.2024196}.

\bibitem[Ananth(2004)]{ranganathan2004levenberg}
Ranganathan Ananth.
\newblock The levenberg-marquardt algorithm.
\newblock 2004.

\bibitem[Wu(1983)]{wu1983convergence}
CF~Jeff Wu.
\newblock On the convergence properties of the em algorithm.
\newblock \emph{The Annals of statistics}, pages 95--103, 1983.

\bibitem[Tarzanagh et~al.(2015)Tarzanagh, Saeidian, Peyghami, and
  Mesgarani]{tarzanagh2015new}
D~Ataee Tarzanagh, Z~Saeidian, M~Reza Peyghami, and H~Mesgarani.
\newblock A new trust region method for solving least-square transformation of
  system of equalities and inequalities.
\newblock \emph{Optimization Letters}, 9:\penalty0 283--310, 2015.

\bibitem[Porcelli(2013)]{porcelli2013convergence}
Margherita Porcelli.
\newblock On the convergence of an inexact gauss--newton trust-region method
  for nonlinear least-squares problems with simple bounds.
\newblock \emph{Optimization Letters}, 7:\penalty0 447--465, 2013.

\end{thebibliography}

\end{document}